\documentclass[12pt,centertags,oneside]{article}
\usepackage{amsmath,amstext,amsthm,amscd,typearea}
\usepackage{amssymb}
\usepackage{a4wide}
\usepackage[mathscr]{eucal}
\usepackage{mathrsfs}
\usepackage{typearea}
\usepackage{charter}
\usepackage{pdfsync}
\usepackage{url}

\usepackage{mathtools}

\usepackage{xcolor}

\usepackage[a4paper,width=16.2cm,top=3cm,bottom=3cm]{geometry}

\usepackage{bbm, dsfont}

\numberwithin{equation}{section}

\newtheorem{theorem}{Theorem}[section]
\newtheorem{definition}[theorem]{Definition}
\newtheorem{proposition}[theorem]{Proposition}
\newtheorem{corollary}[theorem]{Corollary}
\newtheorem{lemma}[theorem]{Lemma}

\newcommand{\cali}[1]{\mathscr{#1}}

\newcommand{\supp}{{\rm Supp}}

\newcommand{\vol}{\mathop{\mathrm{vol}}}

\newcommand{\ddc}{dd^c}
\newcommand{\dc}{d^c}

\newcommand{\Reg}{\text{\normalfont Reg}}
\newcommand{\Sing}{\text{\normalfont Sing}}
\newcommand{\defeq}{\vcentcolon=}

\newcommand{\PSH}{{\rm PSH}}

\newcommand{\C}{\mathbb{C}}

\newcommand{\N}{\mathbb{N}}

\newcommand{\R}{\mathbb{R}}

\renewcommand\P{\mathbb{P}}

\title{\bf Volumes of components of Lelong upper level sets II}
\providecommand{\keywords}[1]{\textbf{\textit{Keywords:}} #1}
\providecommand{\subject}[1]{\textbf{\textit{Mathematics Subject Classification 2010:}} #1}

\author{Shuang Su and Duc-Viet Vu}

\newcommand{\Addresses}{{
                \bigskip
                \footnotesize
                \textsc{Duc-Viet Vu, University of Cologne, Division of Mathematics, Department of Mathematics and Computer Science, Weyertal 86-90, 50931, K\"oln,  Germany}
                \noindent
                \par\nopagebreak
                \noindent
                \textit{E-mail address}: \texttt{dvu@uni-koeln.de} 
                \bigskip
                \par\nopagebreak
                \noindent\textsc{Shuang Su, University of Cologne, Division of Mathematics, Department of Mathematics and Computer Science, Weyertal 86-90, 50931, K\"oln,  Germany}
                \noindent
                \par\nopagebreak
                \noindent
                \textit{E-mail address}: \texttt{ssu1@uni-koeln.de}    
}}

\date{\today}
\begin{document}
\maketitle
\begin{abstract} 
    Let $X$ be a compact K\"ahler manifold of dimension $n$, and let $T$ be a closed positive $(1,1)$-current in a nef cohomology class on $X$. We establish an optimal upper bound for the volume of components of Lelong upper level sets of $T$ in terms of cohomology classes of non-pluripolar self-products of $T$.
\end{abstract}
\noindent
, {density current}, {Lelong number}.
\\

\noindent
\subject{32U15}, {32Q15}.

\section{Introduction}

    The aim of this paper is to investigate the singularities of closed positive currents on compact K\"ahler manifolds.  Let $X$ be a compact K\"ahler manifold of dimension $n$, and let $T$ be a closed positive $(1,1)$-current on $X$. We are interested in understanding the set of points where $T$ has a strictly positive Lelong numbers. By the celebrated upper semi-continuity of Lelong numbers by Siu \cite{Siu}, we know that this set is a countable union of proper analytic subsets on $X$.  Our goal is to estimate the size of this upper level set. We  note that by \cite{Vigny-LelongSkoda}, for every closed positive current $R$ on $X$, there always exists a closed positive $(1,1)$-current whose Lelong numbers coincide with those of $R$. Thus, to study questions about Lelong numbers in this paper, it suffices to consider only currents of bi-degree $(1,1)$. To delve into details, let us first introduce some necessary notions. 

    Let $\omega$ be a fixed smooth K\"ahler form on $X$. We equip $X$ with the Riemannian metric induced by $\omega$. Let $S$ be a closed positive $(p,p)$-current for some $0 \le p \le n$. We define the mass $\|S\|$ of $S$ to be equal to $\int_{X} S \wedge \omega^{n-p}$. For an analytic set $V$ of pure dimension $l$ in $X$, we recall that  
    \[
        \vol(V) =  
        \frac{1}{l !} \int_{\Reg V} \omega^{l},
    \]
    where $\Reg V$ is the regular locus of $V$. Given a closed positive $(p,p)$-current $R$, we denote by $\{R\}  \in H^{p,p}(X,\R)$ the cohomology class of $R$. We say $\alpha \in H^{p,p}(X,\R)$ is \emph{pseudoeffective} if $\alpha =\{R\}$ for some closed positive $(p,p)$-current $R$. Let $\alpha,\beta \in H^{p,p}(X,\R)$. We say $\alpha \geq \beta$ if $\alpha-\beta$ is a pseudoeffective class. 
    
    For every $x \in X$, let $\nu(T,x)$ denote the Lelong number of $T$ at $x$.  For every irreducible analytic subset $V$ in $X$, we recall that the generic Lelong number $\nu(T,V)$ of $T$ along $V$ is defined as $ \inf_{x \in V} \{\nu(T,x)\}$. The Lelong number is a notion measuring the singularities of $T$. We refer to \cite{Demailly_ag} for the basics of Lelong numbers. For every constant $c>0$,  we denote $E_{c}(T) \defeq \{x \in X | \nu(T,x) \ge c\}$ and $E_{+} \defeq \{x \in X | \nu(T,x)>0\}$. By \cite{Siu},  $E_{c}(T)$ is a proper analytic subset in $X$, and $E_{+}(T)= \cup_{m \in \N^*} E_{1/m}(T)$ is a countable union of analytic sets.

    Let $W$ be an irreducible analytic subset of dimension $m$ in $X$. We denote by $E_{+}^{W}(T) \defeq \{x \in W | \nu(T,x) > \nu(T,W)\}$ the Lelong upper level set of $T$ on $W$, which is also a countable union of proper analytic subsets in $W$. 
    Let $V \subset E_{+}^{W}(T)$ be an irreducible analytic set. We say that $V$ is \emph{maximal} if there is no irreducible analytic subset $V'$ of  $E_{+}^{W}(T)$  such that $V$ is a proper subset of $V'$.  We call $V$ \emph{a component of the Lelong upper level set of $T$ along $W$}, and let $\cali{V}_{T,W}$ be the set of such components $V$. Observe that $\cali{V}_{T,W}$ has at most countably many elements. For $0 \leq l \leq m$, we denote by $\cali{V}_{l,T,W}$ the set of $V \in \cali{V}_{T,W}$ such that $\dim V=l$.

    Write $T= \ddc u$ locally, where $u$ is a plurisubharmonic (psh in short) function. We define $T|_{\Reg W}$ to be $\ddc (u|_{\Reg W})$ if $u \not  \equiv -\infty$ on $\Reg W$, and $T|_{\Reg W}:= 0$ otherwise. One sees that this definition is independent of the choice of $u$. Thus, $T|_{\Reg W}$ is a current on $\Reg W$. Here is our main result in the paper.

    \begin{theorem}\label{the-kahlerself-intersec} 
        Let $\alpha$ be a nef $(1,1)$-class and let $W$ be an irreducible analytic subset of dimension $m$ in $X$. Let $T$ be a closed positive current in $\alpha$ such that $\nu(T,W)=0$. Let $1\le m' \le m$ be an integer. Then, we have  
        \begin{align}\label{ine-obstructionselfinterkahler}
            \sum_{V \in \cali{V}_{m-m',T,W}} \nu(T,V)^{m'} \vol(V)  \le 
            \frac{1}{(m-m')!}\int_{\Reg W} \big(\alpha^{m'}   -\langle (T|_{\Reg W})^{m'} \rangle \big) \wedge \omega^{m-m'},
        \end{align}
        where in the integral, we identify $\alpha$ with a smooth closed form in $\alpha$. 
    \end{theorem}

    We have some comments on (\ref{ine-obstructionselfinterkahler}). To see why the term $$I:= \int_{\Reg W} \big(\alpha^{m'}   -\langle (T|_{\Reg W})^{m'} \rangle \big) \wedge \omega^{m-m'}$$ is non-negative, one can consider the case where $W$ is smooth. Then, by a monotonicity of non-pluripolar products (see Theorem \ref{th-monoticity} below), the cohomology class $(\alpha|_W)^{m'}   -\{\langle (T|_{W})^{m'}\rangle\}$ is pseudoeffective. Hence the integral in the right-hand side of (\ref{ine-obstructionselfinterkahler}) is  non-negative. In the general case where $W$ is singular, one can use either a desingularisation of $W$ or interpret $I$ as the mass of some non-pluripolar product relative to $[W]$ (the current of integration along $W$); see Lemma \ref{nonpluripolar_restriction} below.
    We underline however that in order to prove Theorem \ref{the-kahlerself-intersec}, it is not possible to use desingularisation of $W$ to reduce to the case where $W$ is smooth. The reason is that in the process of desingularisation, one has to blow up submanifolds of $W$ which in general could be some of components of the Lelong upper level sets of $T$ on $W$.

    In \cite{Do-Vu-lelong}, a less precise upper bound of volume of components of the Lelong upper level set was  given in terms of the volume of $W$ and the mass of $T$; see also Theorem \ref{th-khisolelongWduong} for a more general statement. The estimate (\ref{ine-obstructionselfinterkahler}) is indeed optimal. To see this, consider $n=m=m'=1$ and $T$ is a sum of Dirac masses. In this case, (\ref{ine-obstructionselfinterkahler}) becomes an equality.     

    If we consider $W=X$, then the generic Lelong number of $T$ along $W$ is zero. Thus, by Theorem \ref{the-kahlerself-intersec}, we have the following result.

    \begin{corollary}    \label{cor-W=X}
        Let $\alpha$ be a nef $(1,1)$-class, and $T$ be a closed positive current in $\alpha$. For $0 \le l \le n$, let $\cali{V}_{l,T}$ be the set of $V \in \cali{V}_{T,X}$ such that $\dim V = l$.
        Let $1\le m' \le n$ be an integer. Then, we have  
        \begin{align*}
            \sum_{V \in \cali{V}_{n-m',T}} \nu(T,V)^{m'} \vol(V)  \le  \frac{1}{(n-m')!}\int_{X} \big{(}\alpha^{m'}  - \{\langle T^{m'} \rangle \} \big{)} \wedge \omega^{n-m'}.
        \end{align*}
    \end{corollary}

    Corollary \ref{cor-W=X} generalizes \cite[Corollary 7.6]{Demailly_regula_11current} by Demailly in which it was assumed additionally that the components of the upper Lelong level set of $T$ are only of dimension $0$ (hence the cohomology class of $T$ is necessarily nef, see \cite[Lemma 6.3]{Demailly_regula_11current}). The feature of Corollary \ref{cor-W=X} is that it holds for any current in a nef class.  

    This paper refines and substitutes \cite{Vu_lelong-bignef-quantitative}. The proof of Theorem \ref{the-kahlerself-intersec} requires both the theory of density currents in \cite{Dinh_Sibony_density} and relative non-pluripolar products in \cite{Viet-generalized-nonpluri} (see also \cite{BT_fine_87,BEGZ}). One of the keys is Theorem \ref{le-truonghom=1lelong} below following from a general comparison of Lelong numbers for density currents (see Corollary \ref{cor-sosanhlelong}).

   The organization of this paper is as follows. In Section \ref{sec-lelong-density}, we recall basic properties of density currents from \cite{Dinh_Sibony_density}. In Section \ref{sec-nonpluripolar}, we discuss the connection between the non-pluripolar product and density currents, and prove Theorem \ref{the-kahlerself-intersec}.\\

    \noindent
    \textbf{Acknowledgments.} This research is partially funded by the Deutsche Forschungsgemeinschaft (DFG, German Research Foundation)-Projektnummer 500055552 and by the ANR-DFG grant QuaSiDy, grant no ANR-21-CE40-0016.

\section{Density currents} \label{sec-lelong-density}

    We first recall some basic properties of density currents introduced by Dinh-Sibony in \cite{Dinh_Sibony_density}. 

    Let $X$ be a complex K\"ahler manifold of dimension $n$, and $V$ a smooth complex submanifold of $X$ of dimension $l.$  Let $T$ be a closed positive $(p,p)$-current on $X,$ where $0 \le p \le n.$  Denote by $\pi: E\to V$ the normal bundle of $V$ in $X$ and $\overline E:= \P(E \oplus \C)$ the projective compactification of $E.$ By abuse of notation, we also use $\pi$ to denote the natural projection from $\overline E$ to $V$. 

    Let $U$ be an open subset of $X$ with $U \cap V \not = \varnothing.$  Let $\tau$ be  a smooth diffeomorphism  from $U$ to an open neighborhood of $V\cap U$ in $E$ which is identity on $V\cap U$ such that  the restriction of its differential $d\tau$ to $E|_{V \cap U}$ is identity.  Such a map is called \emph{an admissible map}.  Note that in \cite{Dinh_Sibony_density}, to define an admissible map,  it is required furthermore that $d\tau$ is $\C$-linear at every point of $V$. This difference doesn't affect what follows.  When $U$ is a small enough tubular neighborhood of $V,$ there always exists an admissible map $\tau$ by \cite[Lemma 4.2]{Dinh_Sibony_density}. In general, $\tau$ is not holomorphic.  When $U$ is a small enough local chart, we can choose a holomorphic admissible map by using suitable holomorphic coordinates on $U$.   For $\lambda \in \C^*,$ let $A_\lambda: E \to E$ be the multiplication by $\lambda$ on fibers of $E.$  Here is the first fundamental result for density currents.

    \begin{theorem} \label{th-dieukienHVconictangenetcurrent} (\cite[Theorem 4.6]{Dinh_Sibony_density})
        Let $\tau$ be an admissible map defined on a tubular neighborhood of $V$. Then,  the family $(A_\lambda)_* \tau_* T$ is of mass uniformly bounded in $\lambda$ on compact subsets in $E$, and  if $S$ is a limit current of the last family as $\lambda \to \infty$, then  $S$ is a current on $E$ which can be extended trivially through $\overline E \backslash E$ to be a closed positive current on $\overline E$  such that the cohomology class $\{S\}$ of $S$ in $\overline E$ is independent of the choice of $S$, and $\{S\}|_V= \{T\}|_V$,  and $\|S\| \le C \|T\|$ for some constant $C$ independent of $S$ and $T$, where $\{S\}|_V$ denotes the pull-back of $\{S\}$ under the natural inclusion map from $V$ to $\overline{E}$.  
    \end{theorem}

    We call $S$  \emph{a tangent current to $T$ along $V$}, and its cohomology class is called \emph{the total tangent class of $T$ along $V$} and is denoted by $\kappa^V(T)$. Tangent currents are not unique in general. By \cite[Theorem 4.6]{Dinh_Sibony_density} again,  if 
    $$S=\lim_{k\to \infty} (A_{\lambda_k})_* \tau_* T$$ for some sequence $(\lambda_k)_k$ converging to $\infty$, then  for every open subset $U'$ of $X$ and  every admissible map $\tau': U' \to E$ , we also have  
    $$S=\lim_{k\to \infty} (A_{\lambda_k})_* \tau'_* T.$$
    This is equivalent to saying that tangent currents are independent of the choice of the admissible map $\tau$.

    \begin{definition} (\cite[Definition 3.1]{Dinh_Sibony_density}) 
        Let $F$ be a complex manifold and $\pi_F: F \to V$ a holomorphic submersion. Let $S$ be  a positive current of bi-degree $(p,p)$  on $F$. \emph{The h-dimension} of $S$ with respect to $\pi_F$ is the biggest integer $q$ such that $S \wedge \pi_F^* \theta^q \not =0$ for some Hermitian metric $\theta$ on $V$.  
    \end{definition}

    By a bi-degree reason, the h-dimension of $S$ is in $[\max\{l- p,0\}, \min\{\dim F -p,l\}]$. 
    We have the following description of currents with minimal h-dimension. 

    \begin{lemma} \label{le-minimalhdimension} (\cite[Lemma 3.4]{Dinh_Sibony_density})
        Let $\pi_F: F \to V$ be a submersion. Let $S$ be a closed positive current of bi-degree $(p,p)$ on $F$ of h-dimension $(l -p)$ with respect to $\pi_F$. Then $S= \pi^* S'$ for some closed positive current $S'$ on $V$. 
    \end{lemma} 

    By \cite[Lemma 3.8]{Dinh_Sibony_density}, the h-dimensions of tangent currents to $T$ along $V$ are the same and this number is called \emph{the tangential h-dimension of $T$ along $V$}.

    Let $m\ge 2$ be an integer. Let $T_j$ be a closed positive current  of bi-degree $(p_j, p_j)$ for $1 \le j \le m$ on $X$ and let  $T_1 \otimes \cdots \otimes T_m$ be the tensor current of $T_1, \ldots, T_m$ which is a current on $X^m.$  A \emph{density current} associated to $T_1, \ldots,  T_m$ is a tangent current to $\otimes_{j=1}^m T_j$ along the diagonal $\Delta_m$ of $X^m.$ Let $\pi_m: E_m \to \Delta$ be the normal bundle of $\Delta_m$ in $X^m$. Denote by $[V]$ the current of integration along $V$.  When $m=2$ and $T_2 =[V]$, the density currents of $T_1$ and  $T_2$ are naturally identified with the  tangent currents to $T_1$ along $V$ (see \cite[Lemma 2.3]{Vu_density-nonkahler}).

    The unique cohomology class of density currents associated to $T_1,\ldots,T_m$ is called \emph{the total density class of $T_1, \ldots, T_m$}. We denote the last class by $\kappa(T_1,\ldots, T_m)$. The tangential h-dimension of  $T_1 \otimes \cdots \otimes T_m$ along $\Delta_m$ is called \emph{the density h-dimension} of  $T_1, \ldots T_m$. 

    \begin{lemma} \label{le-classDSproduct} (\cite[Section 5]{Dinh_Sibony_density})  
        Let $T_j$ be a closed positive current of bi-degree $(p_j,p_j)$ on $X$ for $1 \le j \le m$ such that $\sum_{j=1}^m p_j \le n$. Assume that the density h-dimension of $T_1, \ldots, T_m$ is minimal, \emph{i.e}, equals to $n- \sum_{j=1}^m p_j$. Then the total density class of  $T_1, \ldots, T_m$ is equal to $\pi_{m}^*(\wedge_{j=1}^m\{T_j\})$. 
    \end{lemma}


    Let $h_{\overline E}$ be the Chern class of the dual of the tautological line bundle of $\overline E$. By \cite[Page 535]{Dinh_Sibony_density}, we have  
    \begin{align} \label{eq_kappaVT}
        \kappa^V(T)=\sum_{j=\max \{0, l-p\}}^{\min\{l, n-p-1\}} \pi^*\big( \kappa^V_j(T)\big) \wedge h_{\overline E}^{p-(l-j)},
    \end{align} 
    where $\pi: \overline E \to V$ is the natural projection and $\kappa^V_j(T) \in H^{2l-2j}(V, \R)$. The tangential $h$-dimension of $T$ along $V$ is exactly equal to the maximal $j$ such that $\kappa^V_j(T)  \not = 0$, and it was known that the class $\kappa^V_j(T)$ is pseudoeffective (\cite[Lemma 3.15]{Dinh_Sibony_density}).

    \begin{theorem} \label{th-sosanhVV'densityDS} (\cite[Proposition 4.13]{Dinh_Sibony_density})    Let  $V'$ be a submanifold of $V$ and let $T$ be a closed positive current on $X$. Let       $T_\infty$ be a tangent current to $T$ along $V$. Let $s$ be the tangential h-dimension      of $T_\infty$ along $V'$. Then, the tangential h-dimension of $T$ along $V'$ is at most      $s$, and we have  
            $$\kappa^{V'}_s(T) \le \kappa^{V'}_s(T_\infty).$$
    \end{theorem}

    As a consequence, we obtain the following result generalizing the well-known lower bound of Lelong numbers of intersection of $(1,1)$-currents due to Demailly \cite[Page 169]{Demailly_analyticmethod} in the compact setting. It is probably the first result dealing with comparison of Lelong numbers for intersection of currents of arbitrary bi-degree. 

    \begin{corollary}\label{cor-sosanhlelong} 
        Let $T_j$ be a closed positive current on $X$ for $1 \le j \le m$. Then, for every $x \in X$ and for every density current $S$ associated to $T_1, \ldots, T_{m}$, we have 
        \begin{align}\label{inesosanhSTjlelong}
            \nu(S, x) \ge \nu(T_1, x) \cdots \nu(T_m,x).
        \end{align}
    \end{corollary} 

    \proof  
        Let $x \in X$. Let $\Delta_m$ be the diagonal of $X^m$. Let $\pi: E \to \Delta_m$ be the natural projection from the normal bundle of the diagonal $\Delta_m$ of $X^m$ in $X^m$. We identify $x$ with a point in $\Delta_m$ via the natural identification $\Delta_m$ with $X$.  Put $T:= \otimes_{j=1}^m T_j$ and $V':=\{x\}$. By \cite[Lemma 2.4]{Meo-auto-inter}, we have $\nu(T,x) \ge \nu(T_1,x) \cdots \nu(T_m,x)$. By \cite[Proposition 5.6]{Dinh_Sibony_density}, we have $$\kappa^{V'}_0(S)= \nu(S,x) \delta_{x}, \quad \kappa^{V'}_0(T)= \nu(T,x) \delta_x,$$ where $\delta_{x}$ is the Dirac measure on $x$ (notice here $\dim V' =0$). This combined with Theorem \ref{th-sosanhVV'densityDS} applied to $X^m$, $T:= \otimes_{j=1}^m T_j$, $\Delta_m$ the diagonal of $X^m$ and $V':=\{x\}$  implies $$\nu(S,x)  \ge \nu(T,x) \ge \nu(T_1,x) \cdots \nu(T_m,x).$$ Hence, the desired inequality follows. The proof is finished.  
    \endproof

\section{Relative non-pluripolar products} \label{sec-nonpluripolar}

     We first recall basic facts about the relative non-pluripolar product of currents and discuss its connection with the density current.
 
    Let $X$ be  a compact K\"ahler manifold of dimension $n$. Let $T_1, \ldots, T_m$ be closed positive $(1,1)$-currents on $X$, and $T$ be a closed positive $(p,p)$-current. By \cite{Viet-generalized-nonpluri}, the $T$-relative non-pluripolar product $\langle \wedge_{j=1}^m T_j \dot{\wedge} T\rangle$ is defined  in a way similar to that of  the usual non-pluripolar products  in \cite{BT_fine_87, BEGZ}. For readers' convenience, we explain briefly how to do it. 

    Write $T_j= \ddc u_j+ \theta_j$, where $\theta_j$ is a smooth form and $u_j$ is a $\theta_j$-psh function. Put 
    $$R_k:=\bold{1}_{\cap_{j=1}^m \{u_j >-k\}} \wedge_{j=1}^m (\ddc \max\{u_j,-k\} + \theta_j)\wedge T$$
    for $k \in \N$. By the strong quasi-continuity of bounded psh functions (\cite[Theorems 2.4 and 2.9]{Viet-generalized-nonpluri}), we have 
    $$R_k= \bold{1}_{\cap_{j=1}^m \{u_j >-k\}} \wedge_{j=1}^m (\ddc \max\{u_j,-l\} + \theta_j)\wedge T$$
    for every $l \ge k \ge 1$. One can check that  $R_k$ is positive (see \cite[Lemma 3.2]{Viet-generalized-nonpluri}). 

    As in \cite{BEGZ}, we have that $R_k$ is of mass bounded uniformly in $k$ and $(R_k)_k$ converges to a closed positive current as $k \to \infty$. This limit is denoted by $\langle \wedge_{j=1}^m T_j \dot{\wedge} T\rangle$.  
    We refer to  \cite[Proposition 3.5]{Viet-generalized-nonpluri} for properties of relative non-pluripolar products.

    For every closed positive $(1,1)$-current $P$, we denote by $I_P$ the set of $x \in X$ so that local potentials of $P$ are equal to $-\infty$ at $x$. Note that $I_P$ is  a  complete pluripolar set.  It is clear from the definition that $\langle \wedge_{j=1}^m T_j \dot{\wedge} T\rangle$ has no mass on $\cup_{j=1}^m I_{T_j}$. Furthermore, for every locally complete pluripolar set $A$ (\emph{i.e,} $A$ is locally equal to $\{\psi= -\infty\}$ for some psh function $\psi$), if $T$ has no mass on $A$, then so does $\langle \wedge_{j=1}^m T_j \dot{\wedge} T\rangle$.  The following is deduced from \cite[Proposition 3.5]{Viet-generalized-nonpluri}.

    \begin{proposition}\label{pro-sublinearnonpluripolar} 
        \begin{itemize}
            \item[(i)]For $R:= \langle \wedge_{j=l+1}^m T_j \dot{\wedge} T \rangle$, we have $\langle \wedge_{j=1}^m T_j \dot{\wedge} T \rangle = \langle \wedge_{j=1}^l T_j \dot{\wedge} R \rangle$.

            \item[(ii)]For every complete pluripolar set $A$, we have 
            $$\bold{1}_{X \backslash A}\langle  T_1 \wedge T_2 \wedge \cdots \wedge T_m \dot{\wedge} T\rangle= \big\langle  T_1 \wedge T_2 \wedge \cdots \wedge T_m \dot{\wedge}(\bold{1}_{X \backslash A} T)\big\rangle.$$
            In particular, the equality
            $$\langle \wedge_{j=1}^m T_j \dot{\wedge} T \rangle = \langle \wedge_{j=1}^m T_j \dot{\wedge} T' \rangle$$
            holds, where $T':= \bold{1}_{X \backslash \cup_{j=1}^m I_{T_j}} T$.
        \end{itemize}

    \end{proposition}

    Let $V$ be an irreducible analytic set in $X$. For the case $T=[V]$, we have the following lemma.

    \begin{lemma}\label{nonpluripolar_restriction} (\cite[Lemma 2.3]{Vu-derivative})
        Let $T_1, \dots , T_m$ be closed positive $(1,1)$-currents on $X$. Then the following properties hold:
        \begin{itemize}
            \item[(i)] If $V$ is contained in $\cup_{j=1}^{m}I_{T_j}$, then $\langle T_1 \wedge \dotsi \wedge T_m \dot{\wedge} [V] \rangle=0$ and there is $1 \leq j_0 \leq m$ so that $V \subset I_{T_{j_0}}$.

            \item[(ii)] If $V$ is not contained in $\cup_{j=1}^{m}I_{T_j}$, then 
            \[
                \langle \wedge_{j=1}^{m}  T_j \dot{\wedge} [V] \rangle = i_{*} \langle T_{1,V} \wedge \dotsi \wedge T_{m,V} \rangle,  
            \]
            where $i \colon \Reg(V) \rightarrow X$ is the natural inclusion, and $T_{j,V} \defeq d d^{c}(u_{j}|_{\Reg(V)})$ if $dd^{c}u_j=T_j$ locally.
        \end{itemize}
    \end{lemma}

    Let $T, T'$ be closed positive $(1,1)$-currents in the same cohomology class on $X$. $T'$ is said to be less singular than $T$ if for every local potentials $u$ of $T$ and $u'$ of $T'$, $u \leq u' + O(1)$. Here is a crucial property of relative non-pluripolar products.

    \begin{theorem} \label{th-monoticity} (\cite[Theorem 1.1]{Viet-generalized-nonpluri})   
        Let $T'_j$ be closed positive $(1,1)$-current in the cohomology class of $T_j$ on $X$ such that $T'_j$ is less singular than $T_j$ for $1 \le j \le m$. Then we have 
        $$\{\langle T_1 \wedge \cdots \wedge T_m \dot{\wedge} T \rangle \} \le \{\langle T'_1 \wedge \cdots \wedge T'_m \dot{\wedge} T \rangle\}.$$
    \end{theorem}

    Weaker versions of the above result were proved in  \cite{BEGZ,Lu-Darvas-DiNezza-mono,WittNystrom-mono}. Let $\alpha_1, \dots , \alpha_m$ be big $(1,1)$-classed of $X$. A current $T_{j,\min} \in \alpha_j$ is said to have minimal singularities if it is less singular than any closed positive current in $\alpha_j$.


    By Theorem \ref{th-monoticity}, the class $\{\langle T_{1,\min} \wedge \dotsi \wedge T_{m,\min} \dot{\wedge}T \rangle\}$ is a well-defined pseudoeffective class which is independent of the choice of $T_{j,\min}$. We denote the last class by $\{\langle \alpha_{1} \wedge \dotsi \wedge \alpha_{m} \dot{\wedge} T \rangle \}$. When $T \equiv 1$, we simply write $\langle \alpha_{1} \wedge \dotsi \wedge \alpha_{m} \rangle$ for $\{\langle \alpha_{1} \wedge \dotsi \wedge \alpha_{m} \dot{\wedge} T \rangle \}$. In this case,  the product $\langle \alpha_{1} \wedge \dotsi \wedge \alpha_{m} \rangle $  was introduced in \cite{BEGZ}.  

    Regarding the relation between relative non-pluripolar products and density currents, the following fact was proved in \cite[Theorem 3.5]{Viet-density-nonpluripolar}, see also \cite{VietTuanLucas,Viet_Lucas}.

    \begin{theorem} \label{the-phanbucuarestricteddenstyvanonpluri}  
        Let $R_\infty$ be a density current associated to $T_1, \ldots, T_m,T$. Then we have 
        \begin{align} \label{ine_TjS2}
            \pi_{m+1}^* \langle \wedge_{j=1}^m T_j  \dot{\wedge} T \rangle \le R_\infty,
        \end{align} 
        where $\pi_{m+1}$ is the natural projection from the normal bundle of the diagonal $\Delta$ of $X^{m+1}$ to $\Delta$, and as usual we identified $\Delta$ with $X$. 
    \end{theorem}

    We will need the following to estimate the density h-dimension of currents. 

    \begin{proposition} \label{pro-uocluongmassonsmallsetdensity} (\cite[Proposition 3.6]{Viet-density-nonpluripolar}) 
        Let $T_1,\ldots,T_m,T$ be as above. Let $A$ be a Borel subset of $X$. Assume that for every $J \subset \{1,\ldots, m\}$ and for every density current $R_J$ associated to $(T_j)_{j \in J},T$ with $m_J:=|J|<m$, we have that $R_J$ has no mass on the set $$\pi_{m_J+1}^*(\cap_{j \not \in J} \{x \in A: \nu(T_j, x) >0\}).$$
        Then for every density current $S$ associated to $T_1, \ldots, T_m,T$, the h-dimension of the current $\bold{1}_A S$ is equal to $n-p-m$, where $\bold{1}_A$ denotes the characteristic function of $A$. 
    \end{proposition}

    For every pseudoeffective $(p,p)$-class $\gamma$ on $X$, we put $\|\gamma\|:= \int_X \Theta \wedge \omega^{n-p}$, where $\Theta$ is any closed smooth form  in $\gamma$. This definition is independent of the choice of $\Theta$ and is nonnegative because of the pseudoeffectivity of $\gamma$. 

    \begin{theorem} \label{le-truonghom=1lelong} 
        Let $P$ and $T$ be closed positive currents of bi-degree $(1,1)$ and $(p,p)$ respectively on $X$, where $1 \le p \le n-1$.
        Assume that $T$ has no mass on $I_{P}$. Then, the cohomology class 
        $$\gamma:= \{P\} \wedge \{T\}- \{\langle P \dot{\wedge} T\rangle \}$$
        is pseudoeffective and we have 
        \begin{align}\label{ine-uocluonggammaVPT}
            \| \gamma\| \ge \sum_{V} \nu(P,V)\nu(T,V) 
            n_V!\vol(V),
        \end{align}
        where the sum is taken over every irreducible subset $V$ of dimension
        at least $n-p-1$ in $X$, and $n_V:= \dim V$.   
    \end{theorem}

    \proof  
        Let $\cali{V}$ be the set of irreducible analytic subsets $V$ of dimension at least 
        $n-p-1$ in $X$ such that $\nu(T,V)>0$ and $\nu(P,V)>0$. We note that in (\ref{ine-uocluonggammaVPT}), it is enough to consider  $V \in \cali{V}$. We will see below that $\cali{V}$ has at most countable elements.     

        Observe that if $\nu(P,x)>0$, then $x \in I_P$. Hence, by hypothesis, the trace measure of $T$ has no mass on the set $\{x \in X: \nu(P,x)>0\}$. This allows us to apply Proposition \ref{pro-uocluongmassonsmallsetdensity} to $P$ and $T$ to obtain that the density h-dimension of $P$ and $T$ is minimal. Using this and Lemma \ref{le-classDSproduct} gives 
        \begin{align}\label{eq-densityclassPT}
            \kappa(P,T)= \pi^* (\{P\} \wedge \{T\}),
        \end{align}
        where $\pi$ is the natural projection from the normal bundle of the diagonal $\Delta$ of $X^2$ to $\Delta$.  

        Let $S$ be a density current associated to $P$ and $T$. Since the h-dimension of $S$ is minimal, using Lemma \ref{le-minimalhdimension}, we get that there exists a current $S'$ on $X$ such that $S= \pi^* S'$ (recall $\Delta$ is identified with $X$). Since the relative non-pluripolar product is dominated by density currents (Theorem \ref{the-phanbucuarestricteddenstyvanonpluri}), the current $S' - \langle P \dot{\wedge} T \rangle$ is closed and positive. Moreover, by (\ref{eq-densityclassPT}), the cohomology class of the last current is equal to $\gamma$. It follows that $\gamma$ is pseudoeffective. 

        It remains to prove (\ref{ine-uocluonggammaVPT}).  Let $V \in \cali{V}$. By definition, the generic Lelong number of $T$ along $V$ is positive. Since $T$ is of bi-degree $(p,p)$, the dimension of $V$ must be at most $n-p$. Hence, we have two possibilities: either $\dim V= n-p-1$ or $\dim V= n-p$. The latter case cannot happen. If $T$ has no mass on $V$, then $\nu(T,V)=0$, which leads to a contradiction. If $T$ has mass on $V$, which is contained in $I_P$ (for $\nu(P,V)>0$), then this contradicts the hypothesis that $T$ has no mass on $I_P$. 
  
        Let $V \in \cali{V}$. Since the Lelong numbers are preserved by submersion maps (\cite[Proposition 2.3]{Meo-auto-inter}), by applying Corollary \ref{cor-sosanhlelong} to $P,T$ and generic $x \in V$, we obtain 
        $$\nu(S',V)=\nu(S,V) \ge \nu(P,V) \nu(T,V).$$
        This combined with the fact that $\dim V= n-p-1$ implies  $S' \ge \nu(P,V) \nu(T,V) \, [V]$.
        We deduce that 
        \begin{align*}
            S' &\geq \langle P \dot{\wedge} T \rangle + \bold{1}_{I_P}S' \\
            &\geq  \langle P \dot{\wedge} T \rangle + \sum_{V \in\cali{V}}\nu(P,V) \nu(T,V) \, [V].
        \end{align*}
        The second inequality comes from Siu's decomposition theorem (\cite[2.18]{Demailly_analyticmethod}), and this also shows that $\cali{V}$ has at most countable elements. The desired assertion follows and the proof is finished.
    \endproof

    We now prove Theorem \ref{the-kahlerself-intersec}.

    \begin{proof}[Proof of Theorem \ref{the-kahlerself-intersec}]
        It suffices to consider the case where $\alpha$ is K\"ahler by using $\alpha + \epsilon\{\omega\}, T+ \epsilon \omega$ instead of $\alpha, T$ and letting $\epsilon \to 0$. Hence we assume from now on that $\alpha$ is K\"ahler. By abuse of notation, we also denote by $\alpha$ a smooth K\"ahler form in $\alpha$.  By Lemma \ref{nonpluripolar_restriction}, the right hand-hand side of (\ref{ine-obstructionselfinterkahler}) can be written as 

        \begin{align*}
            \frac{1}{(m-m')!}\int_{\Reg W} \big(\alpha^{m'}   -\langle (T|_{\Reg W})^{m'} \rangle \big) \wedge \omega^{m-m'}= \frac{1}{(m-m')!} \big\|  \langle \alpha^{m'} \dot{\wedge} [W]\rangle- \langle T^{m'} \dot{\wedge} [W] \rangle \big \|.
        \end{align*}

        \mbox{}\\
        \textbf{Step 1.} We first focus on the case that $T$ has analytic singularities. Set $S= \langle T^{m'-1} \dot{\wedge}[W] \rangle$. Since $\alpha$ is K\"ahler, by the monotonicity of non-pluripolar product (Theorem \ref{th-monoticity}) and Proposition \ref{pro-sublinearnonpluripolar} (i), we get 
        \begin{align} \label{ine-sosanhchuyenrelative}
            \big\|  \langle \alpha^{m'} \dot{\wedge} [W]\rangle- \langle T^{m'} \dot{\wedge} [W] \rangle \big \| &\ge 
            \big\|  \langle \alpha \wedge T^{m'-1} \dot{\wedge} [W]\rangle- \langle T^{m'} \dot{\wedge} [W] \rangle \big \|\\
            \nonumber
            &= \| \alpha \wedge S - \langle T \dot{\wedge} S \rangle\|
        \end{align}

        We now show that $S$ has no mass on $I_T$. For $m'>1$, this directly follows from the definition of non-pluripolar product. For $m'=1$, the current $S$ is just $[W]$. Since we assume that $T$ has analytic singularities, the polar locus $I_T$ is an analytic subset and it doesn't contain $W$. Hence, $[W]$ also has no mass on $I_T$. Therefore, we can apply Lemma \ref{le-truonghom=1lelong} to  $T,S$, and get 
        \begin{align}
            \| \alpha \wedge \{S\} - \{\langle T \dot{\wedge} S \rangle\}\|
            &\geq (m-m')! \sum_{V \in \cali{V}_{m-m',T,W}} \nu(T,V) \nu(S,V) \vol(V) \label{3.1}
        \end{align}

        Let $V \in \cali{V}_{m-m',T,W}$ and  let $\Sing(I_T \cap W)$ be the singular locus of the analytic set $I_T \cap W$.  Put $X':= X \backslash \Sing(I_T \cap W)$ which is an open subset in $X$. Now, we claim that the classical intersection $T^{m'-1} \wedge [W]$ is well-defined on some open neighborhood $U$ of $V \backslash \Sing({I_{T} \cap W)}$ in $X'$, in the sense given in \cite[Chapter \uppercase\expandafter{\romannumeral3}]{Demailly_ag}. Since $T$ has analytic singularities, the Lelong number $\nu(T,x)$ is strictly positive if and only if $x$ belongs to $I_T$. This implies that $V$ is contained in $I_{T} \cap W$. Let $K_1, \dots, K_s$ be the irreducible components of $I_{T} \cap W$. Since $V$ is a maximal irreducible proper analytic set such that $\nu(T,V)>0$, $V$ must be one of the irreducible components. By rearranging the index, we may assume $V=K_1$. Let $U$ be an open neighborhood of $V \backslash \Sing(I_{T} \cap W)$ such that 
        \[
            U \cap K_{j}= \emptyset, \mbox{ for }j= 2 , \dots , s.
        \]
        Notice that $V \backslash \Sing(I_T \cap W)$is  contained in $\Reg(V)$, and is of dimension $m-m'$. Consequently, for $0 \leq j' \leq m'-1$,
        \begin{align*}
            \mathscr{H}_{2m-2j'+1}(L(T)|_{U} \cap W) &= \mathscr{H}_{2m-2j'+1}(I_{T} \cap W \cap U)\\
            &= \mathscr{H}_{2m-2j'+1}(V \backslash \Sing(I_{T} \cap W))\\
            &=0,
        \end{align*}
        where $L(T)$ is the set of $x \in X$ such that the local potential of $T$ is unbounded on any neighborhood of $x$. This allows us to apply \cite[Chapter \uppercase\expandafter{\romannumeral3}, Theorem 4.5]{Demailly_ag} and get 
        the well-definedness of $T^{m'-1} \wedge [W] $ on $U$. By \cite[Proposition 3.6]{Viet-generalized-nonpluri}, we obtain 
        \begin{align}
            \langle T^{m'-1} \dot{\wedge} [W] \rangle = \bold{1}_{U \backslash I_{T}} T^{m'-1} \wedge [W]. \label{3.2}
        \end{align}
        Actually, the equality also holds on $U \cap I_T$. To show this, we need to check that $T^{m'-1} \wedge [W]$ has no mass on $U \cap I_{T}$. Since $\dim(U \cap I_{T} \cap W)=m-m'$ and $T^{m'-1} \wedge [W]$ is of bi-dimension $(m-m'+1,m-m'+1)$, the current $T^{m'-1} \wedge [W]$ must have no mass on $U \cap I_{T} \cap W$. Also, by the fact $\supp (T^{m'-1} \wedge [W]) \subset W$, the current $T^{m'-1} \wedge [W]$ also has no mass on $(U \cap I_{T}) \backslash W$. Therefore, the equality (\ref{3.2}) extends to $U$. This implies that the Lelong number $\nu (S,V)$ 
        equals $\nu (T^{m'-1} \wedge [W], V \backslash \Sing(I_T \cap W))$ (remember that we consider the current $T^{m'-1} \wedge [W]$ on $U$, and  $V \backslash \Sing(I_T \cap W)$ is an analytic subset in $U$), and we then have
        \begin{align}
            \nu(S,V) &= \nu (T^{m'-1} \wedge [W], V \backslash \Sing(I_T \cap W)) \notag\\
            &\geq \nu(T,V \backslash \Sing(I_T \cap W))^{m'-1} \nu([W], V \backslash \Sing(I_T \cap W)) \notag\\
            &\geq \nu(T,V)^{m'-1}.\label{3.3}
        \end{align}
        By (\ref{ine-sosanhchuyenrelative}), (\ref{3.1})  and  (\ref{3.3}), the desired inequality follows in the case where $T$ has analytic singularities. 
        \\

        \noindent
        \textbf{Step 2.}
        Now, we remove the assumption that $T$ has analytic singularities. First, we write $T = d \dc u + \theta$, where $\theta$ is a closed smooth $(1,1)$- form and $u \in \PSH(X,\theta)$. Demailly's analytic approximation theorem (see \cite[Corollary 14.13]{Demailly_analyticmethod}) allows us to construct a sequence $u_{k}^{D} \in \PSH(X, \theta + \epsilon_{k} \omega)$, where $\epsilon_k$ decreases to $0$, such that 
        \begin{itemize}
            \item[(1)] $u_{k}^{D} \ge u$ and $u_k^D$ converges to $u$ in $L^1$.
            \item[(2)] $u_{k}^{D}$ has analytic singularities.
            \item[(3)] $\nu(T_k,x)$ converges to $\nu(T,x)$ uniformly on $X$, where $T_k = \ddc u^{D}_k + (\theta + \epsilon_k \omega)$.
        \end{itemize}

        By the monotonicity property of non-pluripolar product (Theorem \ref{th-monoticity}), we have 
        \begin{align} \label{ine0xapxiepsilonTk}
            \big\|  \langle \alpha^{m'} \dot{\wedge} [W]\rangle- \langle T^{m'} \dot{\wedge} [W] \rangle \big \| &= \lim_{k \to \infty} \big\|  \langle (\alpha+ \epsilon_k \omega)^{m'} \dot{\wedge} [W]\rangle- \langle (T+ \epsilon_k \omega)^{m'} \dot{\wedge} [W] \rangle \big \|\\
            \nonumber
            & \ge \limsup_{k\to \infty}\big\|  \langle (\alpha+ \epsilon_k \omega)^{m'} \dot{\wedge} [W]\rangle- \langle T_k^{m'} \dot{\wedge} [W] \rangle \big \|
        \end{align}

        For every constant $r > 0$, set $A_{r}:= \{V \in \cali{V}_{m-m',T,W} | \nu(T,V) \geq r\}$. Observe that $A_r$ increases to $\cali{V}_{m-m',T,W}$ as $r \to 0$.  
        Since $\nu(T_k, x)$ converges to $\nu(T, x)$ uniformly and $T_k$ is less singular than $T$, for every fixed $r>0$ we have 
        \[  
            A_r \subset \cali{V}_{m-m',T_k,W} 
        \]
        when $k$ is large enough. By Step 1, we therefore have 
        \begin{align*}
            \big\|  \langle (\alpha+ \epsilon_k \omega)^{m'} \dot{\wedge} [W]\rangle- \langle T_k^{m'} \dot{\wedge} [W] \rangle \big \|
            &\geq (m-m')!\sum_{V \in \cali{V}_{m-m', T_k,W}} \nu(T_k,V)^{m'} \vol(V)\\
            &\geq (m-m')!\sum_{V \in A_r} \nu(T_k,V)^{m'} \vol(V).
        \end{align*}
        Letting $k \to \infty$ and using (\ref{ine0xapxiepsilonTk}) give
        \begin{align*} 
            \big\|  \langle \alpha^{m'} \dot{\wedge} [W]\rangle- \langle T^{m'} \dot{\wedge} [W] \rangle \big \| &
            \geq (m-m')! \limsup_{k \to \infty}  \sum_{V \in A_{r}} \nu(T_k,V)^{m'} \vol(V)\\ 
            &= (m-m')! \sum_{V \in A_{r}} \nu(T,V)^{m'} \vol(V),
        \end{align*}
        for every constant $r>0$. Letting $r \to 0$, we obtain the desired estimate.  
    \end{proof}

    For the general case where $\nu(T,W)>0$. We couldn't directly compare the volume of Lelong upper level sets of $T$ on $W$ and the mass of $\{\langle \alpha^{m'} \dot{\wedge} [W] \rangle \}-\{\langle T^{m'} \dot{\wedge} [W] \rangle\}\}$. In this case, we have the following modified inequality which is stronger than \cite[Theorem 1.1]{Do-Vu-lelong}.

    \begin{theorem} \label{th-khisolelongWduong}
        Let $\alpha$ be a nef $(1,1)$-class. Let $W$ be an irreducible analytic subset in $X$.  Let $T$ be a closed positive current in $\alpha$. Let $1\le m' \le m$ be an integer.  Then we have
        \begin{multline}\label{ine-vTWduong}
            (m-m')! \sum_{V \in \cali{V}_{m-m',T,W}} \big( \nu(T,V) -\nu(T,W) \big)^{m'} \vol(V) \leq \\
            \big{\|} (\alpha+c \{\omega\})^{m'} \wedge \{[W]\} - \{ \langle (T+c \omega)^{m'} \dot{\wedge} [W] \rangle\}  \big{\|}
        \end{multline}
        where $c= c_1 \cdot \nu(T,W)$ and $c_1>0$ is a constant independent of $\alpha,T,W$. In particular, there is a constant $c_2>0$ independent of $\alpha, T,W$ such that 
        \begin{align} \label{ine-vTWduong2}
            \sum_{V \in \cali{V}_{m-m',T,W}} \big( \nu(T,V) -\nu(T,W) \big)^{m'} \vol(V) \le c_2 \vol(W) \|T\|^{m'}.
        \end{align}
    \end{theorem}

    \begin{proof} The desired inequality (\ref{ine-vTWduong2}) follows directly from (\ref{ine-vTWduong}). We prove now (\ref{ine-vTWduong}). For convenience, set $c_3 := \nu(T,W)>0$. The regularization theorem of currents introduced in \cite{Demailly_regula_11current} allows us to cut down the Lelong upper level set $\{x \in X | \nu (T,V) \geq c_3\}$ from $T$. More precisely, by \cite[Theorem 1.1]{Demailly_regula_11current}, there exists a sequence 
    of almost positive closed $(1,1)$-currents $T_{c_3,k}$ in $\alpha$ such that 
    \begin{itemize}
        \item[(1)] $T_{c_3,k} \geq -(c_1 \cdot c_3 + \epsilon_k) \omega$, where $\lim_{k \rightarrow \infty }\epsilon_k =0$ and $c_1>0$ is a constant independent of $\alpha, T$ and $W$. 
        \item[(2)] The global potentials of $T_{c_3,k}$ decreases to the global potential of $T$.
        \item[(3)] $\nu(T_{c_3,k},x) = \max\{ \nu(T,x)-c_3,0\}$.
    \end{itemize}
    Set $\widetilde{T}_{c_3,k}= T_{c_3,k}+ (c_1 \cdot c_3 + \epsilon_k) \omega$, which is a closed positive $(1,1)$-current. By Theorem \ref{th-monoticity}, we have 
    \begin{align*}
        &\big{\|}(\alpha+ c_1 \cdot c_3 \{\omega\})^{m'} \wedge \{[W]\} - \langle (T+ c_1 \cdot c_3 \omega)^{m'} \dot{\wedge} [W] \rangle \big{\|}\\
        &\geq \limsup_{k \rightarrow \infty} \big{\|}\big{(}\alpha+ (c_1 \cdot c_3+\epsilon_{k}) \{\omega\}\big{)}^{m'} \wedge \{[W]\} - \langle \widetilde{T}_{c_3,k}^{m'}\dot{\wedge} [W] \rangle \big{\|}.
    \end{align*}
    Since $\nu(\widetilde{T}_{c_3,k},W)=0$, we can apply Theorem \ref{the-kahlerself-intersec} to the right-hand side of the above inequality and get 
        \begin{align}
        &\big{\|}\big{(}\alpha+ (c_1 \cdot c_3+\epsilon_{k}) \{\omega\}\big{)}^{m'} \wedge \{[W]\} - \langle \widetilde{T}_{c_3,k}^{m'}\dot{\wedge} [W] \rangle \big{\|} \notag \\
        &\geq (m-m')!\sum_{V \in \cali{V}_{m-m',\tilde{T}_{c_3,k},W}} \nu(\widetilde{T}_{c_3,k},V)^{m'} \vol(V), \label{3.4}
    \end{align}
    By the above properties of $T_{c_3,k}$, we have
    \[
       \cali{V}_{m-m',\tilde{T}_{c_3,k},W} = \cali{V}_{m-m',T,W}.
    \]
    Therefore, the right-hand side of  (\ref{3.4}) is equal to
    \begin{align*}
        &(m-m')!\sum_{V \in \cali{V}_{m-m',T,W}} \nu(\widetilde{T}_{c_2,k},V)^{m'} \vol(V)\\
        &= (m-m')!\sum_{V \in \cali{V}_{m-m',T,W}} \big( \nu(T,V) -\nu(T,W) \big)^{m'}\vol(V).
    \end{align*}
    This completes the proof.
    \end{proof}

\bibliography{biblio_family_MA,biblio_Viet_papers}

\begin{thebibliography}{10}

\bibitem{BT_fine_87}
{\sc E.~Bedford and B.~A. Taylor}, {\em Fine topology, \v{S}ilov boundary, and {$(dd^c)^n$}}, J. Funct. Anal., 72 (1987), pp.~225--251.

\bibitem{BEGZ}
{\sc S.~Boucksom, P.~Eyssidieux, V.~Guedj, and A.~Zeriahi}, {\em Monge-{A}mp\`ere equations in big cohomology classes}, Acta Math., 205 (2010), pp.~199--262.

\bibitem{Lu-Darvas-DiNezza-mono}
{\sc T.~Darvas, E.~Di~Nezza, and C.~H. Lu}, {\em Monotonicity of nonpluripolar products and complex {M}onge-{A}mp\`ere equations with prescribed singularity}, Anal. PDE, 11 (2018), pp.~2049--2087.

\bibitem{Demailly_ag}
{\sc J.-P. Demailly}, {\em Complex analytic and differential geometry}.
\newblock \url{http://www.fourier.ujf-grenoble.fr/~demailly}.

\bibitem{Demailly_regula_11current}
\leavevmode\vrule height 2pt depth -1.6pt width 23pt, {\em Regularization of closed positive currents and intersection theory}, J. Algebraic Geom., 1 (1992), pp.~361--409.

\bibitem{Demailly_analyticmethod}
\leavevmode\vrule height 2pt depth -1.6pt width 23pt, {\em Analytic methods in algebraic geometry}, vol.~1 of Surveys of Modern Mathematics, International Press, Somerville, MA; Higher Education Press, Beijing, 2012.

\bibitem{Dinh_Sibony_density}
{\sc T.-C. Dinh and N.~Sibony}, {\em Density of positive closed currents, a theory of non-generic intersections}, J. Algebraic Geom., 27 (2018), pp.~497--551.

\bibitem{Do-Vu-lelong}
{\sc D.~T. Do and D.-V. Vu}, {\em Volume of components of {L}elong upper-level sets}, J. Geom. Anal., 33 (2023), pp.~Paper No. 303, 14.

\bibitem{VietTuanLucas}
{\sc D.~T. Huynh, L.~Kaufmann, and D.-V. Vu}, {\em Intersection of (1,1)-currents and the domain of definition of the {M}onge-{A}mp\`ere operator}, Indiana Univ. Math. J., 72 (2023), pp.~239--261.

\bibitem{Viet_Lucas}
{\sc L.~Kaufmann and D.-V. Vu}, {\em Density and intersection of {$(1,1)$}-currents}, J. Funct. Anal., 277 (2019), pp.~392--417.

\bibitem{Meo-auto-inter}
{\sc M.~Meo}, {\em In\'{e}galit\'{e}s d'auto-intersection pour les courants positifs ferm\'{e}s d\'{e}finis dans les vari\'{e}t\'{e}s projectives}, Ann. Scuola Norm. Sup. Pisa Cl. Sci. (4), 26 (1998), pp.~161--184.

\bibitem{Siu}
{\sc Y.~T. Siu}, {\em Analyticity of sets associated to {L}elong numbers and the extension of closed positive currents}, Invent. Math., 27 (1974), pp.~53--156.

\bibitem{Vigny-LelongSkoda}
{\sc G.~Vigny}, {\em Lelong-{S}koda transform for compact {K}\"{a}hler manifolds and self-intersection inequalities}, J. Geom. Anal., 19 (2009), pp.~433--451.

\bibitem{Vu_density-nonkahler}
{\sc D.-V. Vu}, {\em Densities of currents on non-{K}\"ahler manifolds}.
\newblock \url{https://doi.org/10.1093/imrn/rnz270}.
\newblock Int. Math. Res. Not. IMRN.

\bibitem{Viet-density-nonpluripolar}
\leavevmode\vrule height 2pt depth -1.6pt width 23pt, {\em Density currents and relative non-pluripolar products}, Bull. Lond. Math. Soc., 53 (2021), pp.~548--559.

\bibitem{Vu_lelong-bignef-quantitative}
\leavevmode\vrule height 2pt depth -1.6pt width 23pt, {\em Loss of mass of non-pluripolar products}.
\newblock \url{ arXiv:2101.05483}, 2021.

\bibitem{Viet-generalized-nonpluri}
\leavevmode\vrule height 2pt depth -1.6pt width 23pt, {\em Relative non-pluripolar product of currents}, Ann. Global Anal. Geom., 60 (2021), pp.~269--311.

\bibitem{Vu-derivative}
\leavevmode\vrule height 2pt depth -1.6pt width 23pt, {\em Derivative of volumes of big cohomology classes}.
\newblock \url{arXiv:2307.15909}, 2023.

\bibitem{WittNystrom-mono}
{\sc D.~Witt~Nystr\"{o}m}, {\em Monotonicity of non-pluripolar {M}onge-{A}mp\`ere masses}, Indiana Univ. Math. J., 68 (2019), pp.~579--591.

\end{thebibliography}
\bibliographystyle{siam}

\bigskip

\noindent
\Addresses
\end{document}